\documentclass{amsart}

\newtheorem{theorem}{Theorem}[section]
\newtheorem{lemma}[theorem]{Lemma}
\newtheorem{corollary}[theorem]{Corollary}

\theoremstyle{definition}
\newtheorem{definition}[theorem]{Definition}

\theoremstyle{remark}
\newtheorem{remark}[theorem]{Remark}

\numberwithin{equation}{section}

\begin{document}

\title[] {On the zeros of Confluent Hypergeometric Functions}

\author{Wei-Chuan  Lin }
\address{Department of Mathematics, Fujian Normal University,\ P. R. China}
\email{sxlwc936@fjnu.edu.cn}

\thanks{This work was supported by the National Natural Science Foundation of China(11371225) and the Natural Science Foundation of Fujian Province in China (2011J01006).}

\author{Xu-Dan  Luo  }
\address{Department of Mathematics, Hong Kong University of
Science and Technology, Clear Water Bay, Kowloon, Hong Kong,\ P. R. China}
\email{xluoab@ust.hk}
\thanks{}

\subjclass[2010]{Primary : 33C15, 30D30, 30D35
}

\keywords{Confluent hypergeometric functions, Jensen Formula, Entire function, Zero-sequence}

\begin{abstract}
In this paper, we study the zero sets of the confluent hypergeometric function $_{1}F_{1}(\alpha;\gamma;z):=\sum_{n=0}^{\infty}\frac{(\alpha)_{n}}{n!(\gamma)_{n}}z^{n}$, where $\alpha, \gamma, \gamma-\alpha\not\in \mathbb{Z}_{\leq 0}$, and show that if $\{z_n\}_{n=1}^{\infty}$ is the zero set of $_{1}F_{1}(\alpha;\gamma;z)$ with multiple zeros repeated and modulus
in increasing order, then there exists a constant $M>0$ such that $|z_n|\geq M n$ for all $n\geq 1$.

\end{abstract}

\maketitle



\section{Introduction}



As we know, many physical, engineering, statistical and applied mathematical problems are related to solutions of a type of second-order linear differential equations, like
\begin{equation}
\label{E:1}
\frac{d^{2}y}{dz^{2}}+p(z)\frac{dy}{dz}+q(z)y(z)=0,
\end{equation}
where $p(z)$ and $q(z)$ are given complex-valued functions.

For example, the investigation of the periods of lateral vibration of a flexible non-uniform rope, or chain, or the periods of vibration of a circular disk, leads to an equation of this type. Moreover, the whirling speed of a non-cylindrical shaft, or the period of lateral vibration of a non-cylindrical bar, such as an air-screw blade, can be found, with two-figure accuracy, by the solution of such an equation; and in fact many vibration problems in various branches of physics lead to such equations. To take another illustration, the crippling end-load of a tapered aeroplane strut, whatever law of taper is adopted, could be found if we could solve equations of this type. Other problems of elastic instability lead to equations of this type, and may be brought into prominence in aeronautics by the urgency of saving weight (see \cite{Webb}).

According to the basic knowledge of Complex Differential Equations, we know the singularities of solutions may only come from those of equation (\ref{E:1}) (see \cite{Bateman}, \cite{Ince}, \cite{G. Kristensson}, \cite{Wang and Guo}).

\begin{definition}
If $z=z_{0}$ is a singularity of (\ref{E:1}), then $z=z_{0}$ is a regular singularity of (\ref{E:1}) if and only if
\begin{equation}
(z-z_{0})\cdot p(z), \ \ (z-z_{0})^{2}\cdot q(z)
\end{equation}
are analytic in $\{z: |z-z_{0}|<R\}$, where $R$ is a positive real number.
\end{definition}

\begin{definition}
A class of equations with $n$ regular singularities is called Fuchsian type equations, where $n$ is a positive integer.
\end{definition}

When $n=3$, the Fuchsian type equation
\begin{equation}
\label{E:2}
z(1-z)\frac{d^{2}y}{dz^{2}}+\left[\gamma-\left(\alpha+\beta+1\right)z\right]\frac{dy}{dz}-\alpha\beta y(z)=0
\end{equation}
is called the Gauss hypergeometric equation, which has three regular singularities $\{0,1,\infty\}$.
\begin{equation}
\label{E:3}
_{2}F_{1}\left(\alpha,\beta;\gamma;z\right):=\sum_{n=0}^{\infty}\frac{(\alpha)_{n}(\beta)_{n}}{n!(\gamma)_{n}}z^{n},\  \gamma\not\in \mathbb{Z}_{\leq 0}
\end{equation}
satisfies (\ref{E:2}) and is convergent in $D(0,1):=\{z:|z|<1\}$, where $(\alpha)_{n}:=\alpha(\alpha+1)\cdots(\alpha+n-1)=\frac{\Gamma(\alpha+n)}{\Gamma(\alpha)}$ is called rising factorial of length $n$.

The confluent hypergeometric equation
\begin{equation}
\label{E:4}
z\frac{d^{2}y}{dz^{2}}+(\gamma-z)\frac{dy}{dz}-\alpha y(z)=0.
\end{equation}
is a confluent case of the Gauss hypergeometric equation by forcing $1$ moving to infinity formally, then $z=0$ is still a regular singularity, but $z=\infty$ becomes an irregular singularity. The process is given by the following.

We replace $z$ by $\frac{z}{b}$ in the Gauss hypergeometric equation (\ref{E:2}), let $b=\beta\rightarrow \infty$, then (\ref{E:2}) will be reduced into the confluent hypergeometric equation (\ref{E:4}).

Around the origin,
\begin{equation}
_{1}F_{1}(\alpha;\gamma;z):=\sum_{n=0}^{\infty}\frac{(\alpha)_{n}}{n!(\gamma)_{n}}z^{n}, \ \gamma\not\in \mathbb{Z}_{\leq 0}
\end{equation}
is one solution of (\ref{E:4}), which is entire, called the confluent hypergeometric function (Kummer function).
In particular, when $\alpha=-n=0,-1,-2,\cdots$, $_{1}F_{1}(-n;\gamma;z)$ is a polynomial (see \cite{G. E. Andrews}, \cite{Bateman}, \cite{Wang and Guo}, \cite{Whittaker and Watson}).

In quantum mechanics, we treat the familiar problem of a quantum harmonic oscillator in one dimension as an application of confluent hypergeometric functions (see \cite{J. Negro}, \cite{J. B. Seaborn}). In statistics, $_{1}F_{1}(\alpha,\gamma;z)$ with integral and half integral values of the parameters $\alpha$ and $\gamma$, occurs in the distributions of many important statistics, such as F-statistics and $D^{2}$-statistics (see \cite{Lebedev}, \cite{P. Nath}).

Besides, many problems in mathematical physics can be solved with the help of the location of zeros of confluent hypergeometric functions. If $\alpha, \gamma, \gamma-\alpha\not\in \mathbb{Z}_{\leq 0}$, then $_{1}F_{1}(\alpha;\gamma;z)$ has infinitely many zeros in $\mathbb{C}$. It is natural to explore the distribution of these zeros.
These are many known results about it for real parameters (see \cite{Slater}). If $\alpha,\gamma\in\mathbb{C}$ are fixed and the variable is sufficiently large, it was showed that the large z-zeros of $_{1}F_{1}(\alpha;\gamma;z)$ satisfy
\begin{equation}
\label{E:zero}
z=\pm(2n+\alpha)\pi i+\ln \left(-\frac{\Gamma(\alpha)}{\Gamma(\gamma-\alpha)}(\pm2n\pi i)^{\gamma-2\alpha}\right)
+O(n^{-1}\ln n),
\end{equation}
where $n$ is a large positive integer, and logarithm takes its principle value (see http://dlmf.nist.gov/13.9). However, this result only tells us the information of zeros with large modulus. In this paper, we would like to look for a common property for each zero.

We denote by $\mathcal{Z}(F)$ the set of zeros of $_{1}F_{1}(\alpha;\gamma;z)$.
Whenever we write $\mathcal{Z}(F)=\{z_n\}$, we always assume that multiple zeros of $_{1}F_{1}(\alpha;\gamma;z)$ are repeated in $\{z_n\}$ and the modulus of the zeros are arranged so that
\begin{equation}
|z_{1}|\leq |z_{2}| \leq \cdots \leq |z_{n}|\leq \cdots.
\end{equation}
Furthermore, we define the infimum of positive number $\tau$ for which $\sum |z_{n}|^{-\tau}$ converges as the exponent of convergence of zero-sequence, and write it as $\lambda(F)$.

It was already known (see\cite{Ahmed}, \cite{Buchhole}) that
\begin{equation}
\label{E:zero-1}
\sum^{\infty}_{j=1}z_j^{-2}=\frac{\alpha(\alpha-\gamma)}{\gamma^2(\gamma+1)}
\end{equation}
and
\begin{equation}
\label{E:zero-2}
\sum^{\infty}_{j=1}z_j^{-3}=\frac{\alpha(\alpha-\gamma)(\gamma-2\alpha)}{\gamma^3(\gamma+1)(\gamma+2)}.
\end{equation}

Afterwards, it was shown by S. Ahmed and M. E. Muldoon (see \cite{S. Ahmed and M. E. Muldoon}) that, if $\{z_n\}$ is a sequence of complex numbers in order of nondecreasing absolute value such that
\begin{equation}
\lim_{n\rightarrow \infty}\frac{z_n}{n}=C\not=0,
\end{equation}
and if $\{z_n\}$ satisfies $(\ref{E:zero-1})$, $(\ref{E:zero-2})$ and
\begin{equation}
2\gamma z^2_k\sum^{\infty}_{j=1}z_j^{-1}(z_k-z_j)^{-1}=(\gamma-2\alpha)z_k-\gamma(\gamma+2),  k=1,2,\cdots,
\end{equation}
then $\{z_n\}$ coincides with the sequence of zeros of $_{1}F_{1}(\alpha;\gamma;z)$. Here, $C$ is a complex number.

Investigating the above results, it is natural to pose the problem: For a confluent hypergeometric function $_{1}F_{1}(\alpha;\gamma;z)$, shall we have $\lambda(F)\leq 1$?

The object of this paper is to devote to the following result.
\begin{theorem}
\label{T:1}
Let $_{1}F_{1}(\alpha;\gamma;z)$ be a confluent hypergeometric function,
where $\alpha, \gamma, \gamma-\alpha\not\in \mathbb{Z}_{\leq 0}$. If $\{z_n\}$ is the zero sequence of
$_{1}F_{1}(\alpha;\gamma;z)$, then there exists a positive constant $M$ such that $|z_n|\geq M n$ for all $n\geq 1$.
\end{theorem}

The condition `$\gamma-\alpha\not\in \mathbb{Z}_{\leq 0}$' in Theorem \ref{T:1} is necessary. For example, we set
$\gamma=\alpha$, then  $_{1}F_{1}(\alpha;\gamma;z)=e^z$, which does not have zeros.

\begin{remark}
The above result for small $n$ can be obtained by (\ref{E:zero}) directly.
\end{remark}

As an immediate consequence of the above theorem, we have

\begin{corollary}
Let $_{1}F_{1}(\alpha;\gamma;z)$ be a confluent hypergeometric function,
where $\alpha, \gamma, \gamma-\alpha\not\in \mathbb{Z}_{\leq 0}$. Then $\lambda(F)\leq 1$ .
\end{corollary}

Jensen's formula (see \cite{Boas}, \cite{GO},\cite{Hayman},\cite{Laine}) represents a property of zeros of analytic functions, which says that
\begin{equation}
\log|f(0)|+\sum^{n}_{k=1}\log{\frac{r}{|a_k|}}=\frac{1}{2\pi}\int^{2\pi}_0\log{|f(re^{i\theta})|d \theta}
\end{equation}
if $f$ is analytic in $|z|\leq r$ with $f(0)\not=0$ and $a_1, a_2,\cdots, a_n$ are the zeros of $f$ in $|z|<r$
repeated according to multiplicity.

\section{Proof of Theorem \ref{T:1}}
In order to verify our main result, we need the following lemma.

\begin{lemma}[\cite{GO}, \cite{Hayman}, \cite{HX}, \cite{Laine}]
\label{L:1} Let $\varphi(x)$ be positive-valued in $[a,b]$, and
$\log \varphi(z)$ be integrable, then
\begin{equation}
\frac{1}{b-a}\int_{a}^{b}\log \varphi(x)dx\leq \log
\left(\frac{1}{b-a}\int_{a}^{b}\varphi(x)dx\right).
\end{equation}
\end{lemma}

Next, we shall prove Theorem \ref{T:1} as follows.
\begin{proof}
Let $\{z_n\} \ \ (n=1,2,\cdots)$ be the zero sequence of the confluent hypergeometric function $f(z):=_{1}F_{1}(\alpha;\gamma;z)$. If $r>0$, then
Jensen's formula gives
\begin{equation}
\log|f(0)|+\sum^{n(r)}_{k=1}\log{\frac{r}{|z_k|}}=\frac{1}{2\pi}\int^{2\pi}_0\log{|f(re^{i\theta})|d \theta},
\end{equation}
where $n(r)$ is the number of zeros of $f$ in $|z|<r$.

Note that $f(0)=1$, exponentiate both sides of the above equality, and apply Lemma \ref{L:1}, we have
\begin{equation}
\prod^{n(r)}_{k=1}\frac{r}{|z_k|}\leq \frac{1}{2\pi}\int^{2\pi}_0|f(re^{i\theta})|d \theta,
\end{equation}
If $n>n(r)$, then
\begin{equation}
\prod^{n}_{k=1}\frac{r}{|z_k|}=\prod^{n(r)}_{k=1}\frac{r}{|z_k|}\prod^{n}_{k=n(r)+1}\frac{r}{|z_k|}\leq
\prod^{n(r)}_{k=1}\frac{r}{|z_k|}.
\end{equation}

If $n<n(r)$, then
\begin{equation}
\prod^{n(r)}_{k=1}\frac{r}{|z_k|}=\prod^{n}_{k=1}\frac{r}{|z_k|}\prod^{n(r)}_{k=n+1}\frac{r}{|z_k|}\geq
\prod^{n}_{k=1}\frac{r}{|z_k|}.
\end{equation}

It follows that
\begin{equation}
\prod^{n}_{k=1}\frac{r}{|z_k|}\leq \frac{1}{2\pi}\int^{2\pi}_0|f(re^{i\theta})|d \theta
\end{equation}
for all $r>0$ and $n\geq 1$.

Since the confluent hypergeometric function
\begin{equation}
f(z):=_{1}F_{1}(\alpha;\gamma;z)=\sum_{m=0}^{\infty}\frac{(\alpha)_{m}}{m!(\gamma)_{m}}z^{m}\ (\gamma\neq 0,-1,-2,\cdots),
\end{equation}
then
\begin{equation}
\begin{split}
&\frac{1}{2\pi}\int^{2\pi}_0|f(re^{i\theta})|d \theta\\
&=\frac{1}{2\pi}\int_{0}^{2\pi}\left|\sum_{m=0}^{\infty}\frac{(\alpha)_{m}}{m!(\gamma)_{m}}r^{m}e^{in\theta}\right|d\theta\\
&\leq\frac{1}{2\pi}\int_{0}^{2\pi}\sum_{m=0}^{\infty}\left|\frac{(\alpha)_{m}}{m!(\gamma)_{m}}\right|r^{m}d\theta\\
&=\sum_{m=0}^{\infty}\left|\frac{(\alpha)_{m}}{m!(\gamma)_{m}}\right|r^{m}.
\end{split}
\end{equation}
Next, we distinguish three cases to discuss.

{\noindent\bf Case 1.} If $\Re \alpha< \Re \gamma$, then we can find a smallest $j\in \mathbb{Z}^{+}$ such that $\Re \alpha+j>0$ and $|\alpha+j|<|\gamma+j|$. Denote
\begin{equation}
C(\alpha,\gamma)=1+\max_{1\leq i\leq j}\left|\frac{(\alpha)_{i}}{(\gamma)_{i}}\right|.
\end{equation}
Thus,
\begin{equation}
\left|\frac{(\alpha)_{m}}{(\gamma)_{m}}\right|\leq C(\alpha,\gamma)
\end{equation}
whenever $0\leq m\leq j$.

Since $\left|\frac{\alpha+k}{\gamma+k}\right|<1$ for $k\geq j$,
we have
\begin{equation}
\left|\frac{(\alpha)_{m}}{(\gamma)_{m}}\right|
=\left|\frac{(\alpha)_{j}}{(\gamma)_{j}}\cdot\frac{(\alpha+j)_{m-j}}{(\gamma+j)_{m-j}}\right|
\leq C(\alpha,\gamma)\cdot \left|\frac{(\alpha+j)_{m-j}}{(\gamma+j)_{m-j}}\right|< C(\alpha,\gamma)
\end{equation}
whenever $j<m$. So
\begin{equation}
\left|\frac{(\alpha)_{m}}{(\gamma)_{m}}\right|\leq C(\alpha,\gamma)
\end{equation}
for $m\geq 0$. It implies that
\begin{equation}
\prod^{n}_{k=1}\frac{r}{|z_k|} \leq \left(C(\alpha,\gamma)\sum_{m=0}^{\infty}\frac{r^{m}}{m!}\right)
= C(\alpha,\gamma) e^{r}.
\end{equation}
{\noindent\bf Case 2.} If $\Re \alpha> \Re \gamma$, then we can find a smallest $\widetilde{j}\in \mathbb{Z}^{+}$ such that $\Re \gamma+\widetilde{j}>0$ and $|\alpha+\widetilde{j}|>|\gamma+\widetilde{j}|$. Thus, there exists a positive integer $\beta\geq 2$ such that
$\beta|\gamma+\widetilde{j}|>|\alpha+\widetilde{j}|$. Similarly, we denote
\begin{equation}
\widetilde{C}(\alpha,\gamma)=1+\max_{1\leq i\leq \widetilde{j}}\left|\frac{(\alpha)_{i}}{(\gamma)_{i}}\right|.
\end{equation}
Then
\begin{equation}
\left|\frac{(\alpha)_{m}}{(\gamma)_{m}}\right|\leq \widetilde{C}(\alpha,\gamma)
\end{equation}
whenever $0\leq m\leq \widetilde{j}$.

Note that $\left|\frac{\alpha+k}{\gamma+k}\right|>1$ for $k\geq \widetilde{j}$
and $\left|\frac{\alpha+k}{\gamma+k}\right|$ is decreasing with respect to $k$, owning lower bound $1$.
So we can obtain that
\begin{equation}
\begin{split}
&\left|\frac{(\alpha)_{m}}{(\gamma)_{m}}\right|=\left|\frac{(\alpha)_{\widetilde{j}}}{(\gamma)_{\widetilde{j}}}
\cdot\frac{(\alpha+\widetilde{j})_{m-\widetilde{j}}}{(\gamma+\widetilde{j})_{m-\widetilde{j}}}\right|\\
&\leq \widetilde{C}(\alpha,\gamma)\cdot \left|\frac{(\alpha+\widetilde{j})_{m-\widetilde{j}}}{(\gamma+\widetilde{j})_{m-\widetilde{j}}}\right|\\
&= \widetilde{C}(\alpha,\gamma)\cdot \left|\frac{(\alpha+\widetilde{j})_{m-\widetilde{j}}}{\beta^{m}\cdot(\gamma+\widetilde{j})_{m-\widetilde{j}}}\right|\cdot\beta^{m}\\
&\leq \widetilde{C}(\alpha,\gamma)\cdot\beta^{m}
\end{split}
\end{equation}
whenever $\widetilde{j}<m$. Then
\begin{equation}
\left|\frac{(\alpha)_{m}}{(\gamma)_{m}}\right|\leq \widetilde{C}(\alpha,\gamma)\cdot\beta^{m}
\end{equation}
for $n\geq 0$. Hence,
\begin{equation}
\prod^{n}_{k=1}\frac{r}{|z_k|} \leq \widetilde{C}(\alpha,\gamma)\sum_{m=0}^{\infty}\frac{(\beta r)^{m}}{m!}
=\widetilde{C}(\alpha,\gamma) e^{\beta r}.
\end{equation}

{\noindent\bf Case 3.} If $\Re \alpha= \Re \gamma$, then there exists a smallest $\widehat{j}\in \mathbb{Z}^{+}$ such that $\Re \alpha+\widehat{j}=\Re \gamma+\widehat{j}>0$. Next, we can follow Case 1 and Case 2 for $|\alpha+\widehat{j}|\leq|\gamma+\widehat{j}|$ and $|\alpha+\widehat{j}|>|\gamma+\widehat{j}|$ respectively.

Combining the above three cases, there exists a constant $C\geq 1$ such that
\begin{equation}
\frac{r^n}{|z_n|^n}\leq \prod^{n}_{k=1}\frac{r}{|z_k|} \leq C e^{\beta r}
\end{equation}
for all $n\geq 1$ and $r>0$.\par
Thus, we have
\begin{equation}
|z_n|^n \geq \frac{r^n}{C e^{\beta r}}
\end{equation}
for  all $n\geq 1$ and $r>0$.\par
We set
$$H(r):=\frac{r^n}{C e^{\beta r}},\ \ \ \ r\in(0, \infty).$$
By a simple calculation, we obtain that
\begin{equation}
\max_{r\in(0, \infty)}H(r)=H\left(\frac{n}{\beta}\right)=C^{-1}\left(\frac{n}{e\beta}\right)^n
\end{equation}
It follows that
\begin{equation}
|z_n|  \geq C^{-\frac{1}{n}}\left(\frac{n}{e\beta}\right)\geq \frac{n}{C e \beta}
\end{equation}
for  all $n\geq 1$.

So we can see that there exists a
constant $M>0$ such that $|z_n|\geq M n$ for all $n\geq 1$.

\end{proof}

\section{Further Remark}

In Section $1$, we have mentioned that, if $\{z_n\}$ is a sequence of complex numbers in order of nondecreasing absolute value such that
\begin{equation}
\label{E:5}
\lim_{n\rightarrow \infty}\frac{z_n}{n}=C\not=0,
\end{equation}
and $\{z_n\}$ satisfies $(\ref{E:zero-1})$, $(\ref{E:zero-2})$ and
\begin{equation}
2\gamma z^2_k\sum^{\infty}_{j=1}z_j^{-1}(z_k-z_j)^{-1}=(\gamma-2\alpha)z_k-\gamma(\gamma+2),  k=1,2,\cdots,
\end{equation}
then $\{z_n\}$ coincides with the sequence of zeros of $_{1}F_{1}(\alpha;\gamma;z)$. Here, $C$ is a complex number.

Combining our main theorem in the paper, we conjecture that $(\ref{E:5})$ can be weakened to
\begin{equation}
|z_n|\geq Mn (M>0)
\end{equation}
for all $n\geq 1$.

\medskip
\bibliographystyle{amsplain}

\end{document}